\pdfoutput=1
\RequirePackage{ifpdf}
\ifpdf 
\documentclass[pdftex]{sigma}
\else
\documentclass{sigma}
\fi

\numberwithin{equation}{section}

\newtheorem{theo}{Theorem}[section]

\newtheorem{lem}[theo]{Lemma}
\newtheorem{prop}[theo]{Proposition}
{\theoremstyle{definition}

\newtheorem{defin}[theo]{Definition}
}

\DeclareMathOperator{\tr}{tr}
\DeclareMathOperator{\dom}{dom}
\DeclareMathOperator{\Ker}{Ker}
\begin{document}

\allowdisplaybreaks

\renewcommand{\PaperNumber}{036}

\FirstPageHeading

\ShortArticleName{A Note on Gluing Dirac Type Operators on a~Mirror Quantum Two-Sphere}

\ArticleName{A Note on Gluing Dirac Type Operators\\
on a~Mirror Quantum Two-Sphere}

\Author{Slawomir KLIMEK~$^\dag$ and Matt MCBRIDE~$^\ddag$}

\AuthorNameForHeading{S.~Klimek and M.~McBride}

\Address{$^\dag$~Department of Mathematical Sciences, Indiana University-Purdue University Indianapolis,\\
\hphantom{$^\dag$}~402 N.~Blackford St., Indianapolis, IN 46202, USA}
\EmailD{\href{mailto:sklimek@math.iupui.edu}{sklimek@math.iupui.edu}}

\Address{$^\ddag$~Department of Mathematics, University of Oklahoma, 601 Elm St., Norman, OK 73019, USA}
\EmailD{\href{mailto:mmcbride@math.ou.edu}{mmcbride@math.ou.edu}}

\ArticleDates{Received September 30, 2013, in f\/inal form March 25, 2014; Published online March 29, 2014}

\Abstract{The goal of this paper is to introduce a~class of operators, which we call quantum Dirac type operators on
a~noncommutative sphere, by a~gluing construction from copies of noncommutative disks, subject to an appropriate local
boundary condition.
We show that the resulting operators have compact resolvents, and so they are elliptic operators.}


\Classification{46L99; 47B99; 81R60}

\section{Introduction}

This paper is a~continuation of our study of operators that we are constructing using noncommutative geometry concepts,
and which we named quantum Dirac type operators on noncommutative spaces, see~\cite{CKW,KM1,KM2,KM4,KM5}.
The goal of our previous papers was to provide simple examples of such operators on noncommutative compact manifolds
with boundary and then study non-local, Atiyah--Patodi--Singer type boundary conditions and the corresponding index
problem.
In particular we concentrated our ef\/forts on showing that those operators are invertible or invertible modulo compact
operators, and that the inverses/parametrices are also compact operators, which are the essential properties of
ellipticity.

In this paper we def\/ine an interesting class of operators which we call as quantum Dirac type operators on
a~noncommutative sphere of~\cite{HMS}, by gluing two operators from noncommutative disks, subject to an appropriate
boundary condition.
In doing so we are using the concepts and techniques from our previous papers, even though in this note we employ
a~local boundary condition.
Other constructions of particular Dirac operators on (Podle\'s) quantum spheres are described
in~\cite{BL,DLPS,DS,DDLW}, using dif\/ferent techniques.
It has to be stressed that in this paper we do not single out any one Dolbeault or Dirac operator, but rather we provide
a~functional analytic framework for studying a~large class of such operators, an even more general class than we
introduced in~\cite{KM5}.
We use the terminology ``quantum Dirac type operators" because of many similarities between them and classical Dirac
type operators, including their structure in the Fourier transform, as seen by comparing the formulas~\eqref{classDTO}
and~\eqref{quantDTO}.
It is a~very interesting problem to cast those operators in the spectral triples framework and we plan to address it in
future publications.
The main goal, as in the past papers, is to show that our quantum Dirac type operators have compact resolvents.
We try to be quite self-contained in the presentation so that this paper can be read without referring too much to our
previous work.

In Section~\ref{Section2} we describe how to glue together two d-bar operators on the unit disk to make the d-bar operator on the
two-sphere in classical dif\/ferential geometry.
The purpose here is to understand the boundary condition which needs to be incorporated in the gluing.
This is followed by the computation of the kernel of the d-bar operator and the construction of its parametrix.
The calculations then lead to an alternative way of proving the standard results on Fredholm theory and the resolvent of
the d-bar operator on the sphere.

In Section~\ref{Section3} we def\/ine the quantum two-sphere of~\cite{HMS}, the quantum disk of~\cite{KL}, our class of quantum Dirac
type operators, and the Hilbert spaces playing the role of the $L^2$ spaces of the previous section.
Just like in the commutative case, the boundary condition will come from how the quantum two-sphere is formed.
The calculations of the kernel and the parametrix for the quantum Dirac type operator parallel the classical case and
are quite similar to those of our previous papers.
The proof of the main result on compactness of the parametrices follows.
Finally at the end of the section a~specif\/ic example of the so-called $q$-weights is given that is motivated
by~\cite{KM1}.

\section{Gluing disks classically}\label{Section2}

In this section we describe the classical dif\/ferential geometry needed to create the two-sphere and the d-bar operator
on it via gluing.
We consider the standard, round two-sphere or, equivalently, the Riemann extended plane equipped with the Fubini--Study
metric, i.e.~the constant curvature metric, see~\cite{Ahlfors} for details.
Let
\begin{gather*}
{\mathbb D} = \{z\in{\mathbb C}: |z|\le 1\}
\end{gather*}
be the disk with radius one.
Identify the outside of ${\mathbb D}$ with ${\mathbb D}$ via
\begin{gather*}
z\mapsto\frac{1}{z},
\end{gather*}
and so the Riemann sphere becomes the union of two unit disks modulo the gluing on the boundary.

The Hilbert space of functions on the sphere is then identif\/ied with the direct sum
\begin{gather*}
\mathcal{H}_1 = L^2({\mathbb D}, d\mu) \oplus L^2({\mathbb D}, d\mu),
\end{gather*}
where the measure $d\mu$ is given by the formula
\begin{gather*}
d\mu = \frac{1}{(1+|z|^2)^2}d^2z
\end{gather*}
with $d^2z$ being the Lebesgue measure.
The Hilbert space of $(0,1)$ dif\/ferential forms is
\begin{gather*}
\mathcal{H}_2 = L^2\big({\mathbb D}, d^2z\big) \oplus L^2\big({\mathbb D}, d^2z\big).
\end{gather*}

The geometric d-bar operator takes functions into $(0,1)$ dif\/ferential forms: $f\to (\overline{\partial}f)d\bar z=
(\partial f/\partial\overline{z})d\bar z$.
The d-bar operator $D: \mathcal{H}_1\to\mathcal{H}_2$ on the sphere is def\/ined naturally to be
\begin{gather}
\label{dirac_op}
D(u,v) = (\overline{\partial}u, \overline{\partial}v)
\end{gather}
for $(u,v)\in\mathcal{H}_1$.
The domain of $D$ is given by
\begin{gather}
\label{dom_D}
\dom(D) = \left\{(u,v)\in\mathcal{H}_1: u,v\in H^1({\mathbb D})~\text{and}~u(z)=v(1/z)~\textrm{when}~|z|=1\right\},
\end{gather}
where $H^1({\mathbb D})$ is the f\/irst Sobolev space on the disk.

The key tool that we use to analyze $D$ is the following Fourier decomposition for $u\in L^2({\mathbb D},d^2z)$
\begin{gather*}
u(z) = \sum\limits_{n\ge0} u_n^+(r) e^{in\theta} + \sum\limits_{n\ge1}u_n^-(r)e^{-in\theta},
\end{gather*}
and the same for $v\in L^2({\mathbb D},d\mu)$.
This leads to natural decompositions of the Hilbert spaces
\begin{gather*}
L^2\big({\mathbb D}, d^2z\big) \cong \bigoplus_{n\ge0} L^2\left([0,1],rdr\right)\oplus \bigoplus_{n\ge
1}L^2\left([0,1],rdr\right)
\end{gather*}
and
\begin{gather*}
L^2({\mathbb D},d\mu(z)) \cong \bigoplus_{n\ge0} L^2\left([0,1], d\mu(r)\right) \oplus \bigoplus_{n\ge1} L^2\left([0,1],d\mu(r)\right),
\end{gather*}
where
\begin{gather*}
d\mu(r):= \frac{r}{(1+r^2)^2}dr.
\end{gather*}

The point is that the Fourier decomposition naturally carries forward to the quantum disk case with little
modif\/ications.
It will be useful to write the d-bar operator in the Fourier components
\begin{gather}
\label{classDTO}
\overline{\partial}u(z) = \frac{1}{2}\sum\limits_{n\ge0}\!\left(\left(u_n^+\right)'(r) - \frac{n}{r}u_n^+(r)\right)e^{i(n+1)\theta}
+ \frac{1}{2}\sum\limits_{n\ge1}\!\left(\left(u_n^-\right)'(r) + \frac{n}{r}u_n^-(r)\right)e^{-i(n-1)\theta}.
\end{gather}

Using the Fourier decompositions of both $u$ and $v$, the boundary condition $u(z)=v(1/z)$ can be rewritten as
\begin{gather}
\label{bcft}
u_n^+(1) = v_n^-(1)
\qquad
\textrm{and}
\qquad
u_n^-(1) = v_n^+(1),
\qquad
n\ge 1,
\qquad
\textrm{and}
\qquad
u_0^+(1) = v_0^+(1).
\end{gather}

The above descriptions of the d-bar operator on the sphere and the Fourier expansion allow us to prove the well known
functional analytic properties of $D$.

\begin{theo}\label{compact_para_cont}
Let $D$ be the Dirac operator defined by~\eqref{dirac_op} on domain~\eqref{dom_D}, then $D$ is a~right invertible
operator and it is left invertible modulo compacts.
Moreover $D$ has a~compact parametrix.
\end{theo}

We begin by computing the kernel of $D$.

\begin{prop}
The kernel of the d-bar operator $D$ defined by~\eqref{dirac_op} on domain~\eqref{dom_D} is one-dimensional and consists
of constant functions.
\end{prop}

\begin{proof}
We need to solve the following uncoupled system
\begin{gather*}
\overline{\partial}u = 0,
\qquad
\overline{\partial}v = 0
\end{gather*}
subject to the boundary condition $u(z)=v(1/z)$ when $|z|=1$.
Using the Fourier decomposition in the f\/irst equation leads to
\begin{gather*}
(u_n^+)'(r) - \frac{n}{r} u_n^+(r) = 0
\qquad
\textrm{and}
\qquad
(u_n^-)'(r) + \frac{n}{r} u_n^-(r) = 0.
\end{gather*}
The solution to the f\/irst equation is $u_n^+(r) = c_n^+ r^n$ for constants $c_n^+$.
Moreover $u_n^-(r) = c_n^-r^{-n}$ for constants $c_n^-$.
Similarly we have $v_n^+(r) = d_n^+r^n$ and $v_n^-(r) = d_n^-r^{-n}$ for constants $d_n^+$ and $d_n^-$.
Since the solutions are required to be in $H^1({\mathbb D})$ we must have $c_n^- = d_n^- = 0$ for all $n\ge 1$.
Now we apply the boundary condition~\eqref{bcft}.
For $n\ge1$ we obtain $c_n^+ = d_n^-$ and $d_n^+ = c_n^-$ which means that $c_n^+ = 0$ and $d_n^+ = 0$.
If $n=0$ we get $c_0^+=d_0^+$ is an arbitrary constant.
This tells us that $\Ker D$ is one-dimensional and the proof is complete.
\end{proof}

The next problem is to f\/ind an operator $Q$ so that $DQ(u,v)=(u,v)$.
To write a~formula for the right inverse $Q$ we need the following three integral operators.
Let $\chi(t)=1$ for $t\le1$ and zero otherwise.
Def\/ine
\begin{gather}
T_1^{(n)}f(r) = \int_0^1 2r^n\rho^{n-1}f(\rho)\rho d\rho,
\qquad
T_2^{(n)}f(r) = -\int_0^1 2\chi\left(\frac{r}{\rho}\right)\frac{r^n}{\rho^{n+1}}f(\rho)\rho d\rho,
\nonumber
\\
T_3^{(n)}f(r) = \int_0^1 2\chi\left(\frac{\rho}{r}\right)\frac{\rho^{n-1}}{r^n}f(\rho)\rho d\rho
\label{int_ops_cont}
\end{gather}
and $T_1^{(n)}$, $T_2^{(n)}$, and $T_3^{(n)}: L^2([0,1], rdr)\to L^2([0,1], d\mu(r))$.

\begin{prop}
\label{cont_para}
Let the operator $D$ be given by~\eqref{dirac_op} on domain~\eqref{dom_D}.
The operator $D$ has a~right inverse $Q$, i.e.
$DQ(p,q)=(p,q)$ where $Q(p,q)=({u},{v})$ and
\begin{gather*}
{u} = T_2^{(0)}p_1^+(r) + \sum\limits_{n\ge 1}\left(T_2^{(n)}p_{n+1}^+(r) + T_1^{(n)}q_{n-1}^-(r)\right)e^{in\theta} +
\sum\limits_{n\ge1} T_3^{(n)}p_{n-1}^-(r)e^{-in\theta},
\\
{v} = T_2^{(0)}q_1^+(r) + \sum\limits_{n\ge 1}\left(T_2^{(n)}q_{n+1}^+(r) + T_1^{(n)}p_{n-1}^-(r)\right)e^{in\theta} +
\sum\limits_{n\ge1} T_3^{(n)}q_{n-1}^-(r)e^{-in\theta}.
\end{gather*}
Moreover we have $QD= I - C$, where $C(u,v) = \left(u_0^+(1), u_0^+(1)\right)$.
\end{prop}

\begin{proof}
In order to construct a~parametrix of $D$, we need to solve the equation $D(u,v) = (p,q)$ with $(u,v)$ subject to the
boundary condition.
This is the same as solving the uncoupled system
\begin{gather*}
\overline{\partial}u = p,
\qquad
\overline{\partial}v = q
\end{gather*}
subject to the condition~\eqref{bcft}.
To solve $\overline{\partial}u=p$ we use the Fourier decomposition and the equation reduces to
\begin{gather*}
\frac{1}{2}\left((u_n^+)'(r) - \frac{n}{r}u_n^+(r)\right) = p_{n+1}^+(r)
\qquad
\textrm{and}
\qquad
\frac{1}{2}\left((u_n^-)'(r) + \frac{n}{r}u_n^-(r)\right) = p_{n-1}^-(r).
\end{gather*}
We obtain the following solutions
\begin{gather*}
u_n^+(r) = c_n^+ r^n - 2\int_{r}^{1} \left(\frac{r}{\rho}\right)^n p_{n+1}^+(\rho)d \rho
\qquad
\textrm{and}
\qquad
u_n^-(r) = 2\int_{0}^{r} \left(\frac{\rho}{r}\right)^n p_{n-1}^-(\rho)d \rho
\end{gather*}
for some constants $c_n^+$ to be computed from the boundary condition.
Here the limits of integration are determined so that the functions are non singular at $r=0$ to be in the domain of~$\overline{\partial}$.
The solution for $v_n^+(r)$ and $v_n^-(r)$ are given by identical formulas with the corresponding constants denoted by
$d_n^+$.
If $n\ge1$ then applying $u_n^+(1) = v_n^-(1)$ we get
\begin{gather*}
c_n^+ = u_n^+(1) = v_n^-(1) = 2\int_0^1 \rho^n q_{n-1}^-(\rho)d\rho
\end{gather*}
and similarly applying $v_n^+(1) = u_n^-(1)$ we obtain
\begin{gather*}
d_n^+ = v_n^+(1) = u_n^-(1) = 2\int_0^1 \rho^n p_{n-1}^-(\rho)d\rho.
\end{gather*}
Finally $c_0^+=d_0^+$.
Comparing the above formulas with the integrals~\eqref{int_ops_cont}, we get
\begin{gather*}
u_n^+(r) = T_1^{(n)}q_{n-1}^-(r) + T_2^{(n)}p_{n+1}^+(r),
\qquad
u_n^-(r) = T_3^{(n)}p_{n-1}^-(r),
\\
v_n^+(r) = T_1^{(n)}p_{n-1}^-(r) + T_2^{(n)}q_{n+1}^+(r),
\qquad
v_n^-(r) = T_3^{(n)}q_{n-1}^-(r).
\end{gather*}
It remains to check that $QD(u,v) = (I-C)(u,v)$, which is a~direct but somewhat lengthly computation and will be
omitted.
This f\/inishes the proof.
\end{proof}

An immediate corollary is that since the operator $C$ def\/ined in the above proposition is a~f\/inite rank operator, $D$ is
left invertible modulo compacts.

The next goal is to show that $Q$ is a~compact operator.
We f\/irst investigate the compactness and norms of the integral operators $T_1^{(n)}$, $T_2^{(n)}$, and $T_3^{(n)}$.
This is summarized in the following proposition.

\begin{prop}
\label{hs_norm_cont_T}
The three integral operators given by formula~\eqref{int_ops_cont} are Hilbert--Schmidt ope\-rators.
Moreover the Hilbert--Schmidt norms of these integral operators go to zero as $n\to\infty$.
\end{prop}

\begin{proof}
We f\/irst assume $n\ge1$ for all of the operators and then consider the operator $T_2^{(0)}$ separately.
For all the following estimates we use the fact that $(1+r^2)^{-2} \le 1$ for $0\le r\le 1$.
The Hilbert--Schmidt norm of $T_1^{(n)}$ is bounded by
\begin{gather*}
\|T_1^{(n)}\|_{\textrm{HS}}^2 \le 4\int_0^1\int_0^1r^{2n+1}\rho^{2n-1}d\rho dr =\frac{1}{n(n+1)}.
\end{gather*}
Similarly
\begin{gather*}
\|T_2^{(n)}\|_{\textrm{HS}}^2 \le 4\int_0^1\int_r^1 \left(\frac{r}{\rho}\right)^{2n+1}d\rho dr=4\left(\frac{1}{4n}-
\frac{1}{4n(n+1)}\right) \le \frac{1}{n},
\end{gather*}
and
\begin{gather*}
\|T_3^{(n)}\|_{\textrm{HS}}^2 \le 4\int_0^1\int_0^r \frac{\rho^{2n-1}}{r^{2n-1}}d\rho dr = \frac{1}{n}.
\end{gather*}
Finally we have
\begin{gather*}
\big\|T_2^{(0)}\big\|_{\textrm{HS}}^2 \le 4 \int_0^1 \int_r^1 \frac{r}{\rho}d\rho dr = 4\lim_{t\to0^+}\int_1^t r\ln r dr
<\infty
\end{gather*}
since $r\ln r$ is continuous at $r=0$.
It is clear that the Hilbert--Schmidt norms are all f\/inite and go to zero as $n\to\infty$.
This completes the proof.
\end{proof}

\begin{proof}[Proof of Theorem~\ref{compact_para_cont}] Proposition~\ref{cont_para} shows that $Q$ is the right inverse to $D$.
Moreover the comment directly after the proof of Proposition~\ref{cont_para} shows that $D$ is left invertible modulo
compact operators.
We use the decomposition of $Q$ in Proposition~\ref{cont_para} and the fact that the Hilbert--Schmidt norms of
$T_1^{(n)}$, $T_2^{(n)}$, and $T_3^{(n)}$ go to zero as $n\to\infty$ in Proposition~\ref{hs_norm_cont_T}, to conclude
that $Q$ must be a~compact operator as the norm limit of compact operators, its partial sums.
This is because the operators $T_1^{(n)}$, $T_2^{(n)}$, and $T_3^{(n)}$ act in mutually orthogonal subspaces of the
Hilbert space and so $Q$ is an inf\/inite direct sum of compact operators whose norms are going to zero.
This ends the proof.
\end{proof}

\section{Quantum sphere}\label{Section3}

In noncommutative geometry the quantum spaces are def\/ined to be specif\/ic $C^*$-algebras, which play the role of the
spaces of continuous functions.
With this in mind we def\/ine the quantum disk.

\begin{defin}
Let $\{e_k\}$ be the canonical basis for $\ell^2({\mathbb N})$ and let $U$ be the unilateral shift, i.e.
$Ue_k = e_{k+1}$.
Let $C^*(U)$ be the $C^*$-algebra generated by $U$.
The $C^*$-algebra $C^*(U)$ is called the quantum unit disk~\cite{KL}.
\end{defin}

Below we will need some results about $C^*(U)$ which we will review now.
There is a~short exact sequence
\begin{gather*}
0\longrightarrow \mathcal{K} \longrightarrow C^*(U) \overset{\sigma}{\longrightarrow} C(S^1) \longrightarrow 0,
\end{gather*}
where $ \mathcal{K} $ is the algebra of compact operators in $\ell^2({\mathbb N})$, and it is the commutator ideal of
$C^*(U)$.
The commutative quotient $C^*(U)/\mathcal{K}$ is isomorphic with $C(S^1)$ where the quotient map $\sigma$, also called
the restriction to the boundary map is given by $\sigma(U)= e^{i\theta}\in C(S^1) $ and $\sigma(U^*)=e^{-i\theta}$.
The related $C^*$-algebra morphism $\overline{\sigma}:C^*(U)\mapsto C(S^1)$, def\/ined by the formulas
$\overline{\sigma}(U)=e^{-i\theta}$ and $\overline{\sigma}(U^*) = e^{i\theta}$, is the composition of $\sigma$ with the
$z\to 1/z$ automorphism on the unit circle.

Let $Ke_k = ke_k$ be the label operator.
The operators $(K,U)$ are called noncommutative polar coordinates.

For a~numerical function $f(k)$, the diagonal operator $f(K)$ in $\ell^2({\mathbb N})$ belongs to $C^*(U)$ if\/f
$\lim\limits_{k\to\infty}f(k)$ exists, and in fact we have
\begin{gather*}
\sigma(f(K))=\lim_{k\to\infty}f(k)I\in C(S^1).
\end{gather*}
The same formula is true for $\overline{\sigma}$.

For $x\in C^*(U)$ def\/ine the following formal series of operators
\begin{gather*}
x_\textrm{series} = \sum\limits_{n\ge0} U^n x_n^+(K) + \sum\limits_{n\ge1} x_n^-(K)(U^*)^n,
\end{gather*}
where $x_n^+(k) = \langle e_k, (U^*)^nx e_k\rangle$ and $x_n^-(k) = \langle e_k, xU^n e_k\rangle$.

Similarly to the usual theory of Fourier series, $x_\textrm{series}$ determines $x$ even though in general the series is
not norm convergent.
Other types of convergence results can be obtained along the lines of the usual Fourier analysis.

We have
\begin{gather}
\label{lim_at_infty1}
\sigma(x) = \sum\limits_{n\ge0} e^{in\theta}x_n^+(\infty) + \sum\limits_{n\ge1} x_n^-(\infty)e^{-in\theta}
\end{gather}
and
\begin{gather}
\label{lim_at_infty2}
\overline{\sigma}(y) = \sum\limits_{n\ge0} e^{-in\theta} y_n^+(\infty) + \sum\limits_{n\ge1} y_n^-(\infty)e^{in\theta},
\end{gather}
where $f(\infty):=\lim\limits_{k\to\infty}f(k)$, provided the limit exists.

We now def\/ine the Hilbert space that plays the role of the space of square integrable functions or the square integrable
$(0,1)$ dif\/ferential forms on the disk.
While for the dif\/ferential geometric reason it is important to have appropriate measures in the $L^2$ spaces of
functions and one-forms, those details are not important for the kind of functional analytic properties of the d-bar
operator investigated in this paper, as long as the measures are equivalent.
For simplicity in the quantum case we consider a~general class of Dirac like operators acting in a~f\/ixed Hilbert space.

In the def\/inition of that space we use Fourier series, and the f\/inite volume property of the disk is achieved using
weights.
Let $\{a^{(n)}(k)\}$ be a~sequence of positive numbers such that the sum $s(n):=\sum\limits_{k=0}^\infty
\frac{1}{a^{(n)}(k)}$ exists and such that $s(n)$ goes to zero as $n\to\infty$.
The f\/irst condition ref\/lects in a~sense the f\/inite volume of the disk.
Notice that the weights studied in~\cite{KM1} decrease monotonically in $n$ for each $k$.
In this paper we allow more general weights than in that paper.

For a~formal power series
\begin{gather*}
f = \sum\limits_{n\ge0} U^n f_n^+(K) + \sum\limits_{n\ge1} f_n^-(K)(U^*)^n
\end{gather*}
we def\/ine the following weighted $\ell^2$ norm
\begin{gather*}
\|f\|^2 = \sum\limits_{n\ge 0} \sum\limits_{k=0}^\infty \frac{1}{a^{(n)}(k)}|f_n^+(k)|^2 +
\sum\limits_{n\ge1}\sum\limits_{k=0}^\infty \frac{1}{a^{(n)}(k)}|f_n^-(k)|^2.
\end{gather*}
Let $\mathcal{H}_1$ be the Hilbert space whose elements are power series $f$ such that $\|f\|$ is f\/inite.
We have the following density proposition.

\begin{prop}
If $x\in C^*(U)$, then $x_{\rm series}$ converges in $\mathcal{H}_1$.
Moreover the map $C^*(U)\ni x\to x_{\rm series}\in \mathcal{H}_1$ is continuous, one-to-one, and the image is dense
in $\mathcal{H}_1$.
\end{prop}

\begin{proof}
For $x\in C^*(U)$ we need to estimate the norm of $x_\textrm{series}$.
Notice f\/irst that if $x$ is a~f\/inite sum
\begin{gather*}
x = \sum\limits_{n=0}^N U^n x_n^+(K) + \sum\limits_{n=1}^N x_n^-(K)(U^*)^n
\end{gather*}
then $x_{\textrm{series}} = x$.
Since such $x$'s are dense in $C^*(U)$, it suf\/f\/ices to estimate the norm of the f\/inite sums.
For brevity, we assume that $x$ has only the $U^*$ terms since the same argument will work for the $U$ terms as well.
We have
\begin{gather*}
\|x_{\textrm{series}}\|_{\mathcal{H}_1}^2 = \sum\limits_{k=0}^\infty \sum\limits_{n=1}^N
\frac{1}{a^{(n)}(k)}|x_n^-(k)|^2 = \tr\left(\sum\limits_{n=1}^N \frac{1}{a^{(n)}(K)}|g_n(K)|^2\right)
\\
\phantom{\|x_{\textrm{series}}\|_{\mathcal{H}_1}^2}
= \tr\left(\sum\limits_{n=1}^N \frac{1}{a^{(n)}(K)}g_n(K)\left(U^*\right)^nU^n\overline{g_n(K)}\right)
\\
\phantom{\|x_{\textrm{series}}\|_{\mathcal{H}_1}^2}
=\tr\left(\sum\limits_{n=1}^N \frac{1}{a^{(n)}(K)}g_n(K)\left(U^*\right)^n\sum\limits_{l=1}^NU^l\overline{g_l(K)}\right)
\\
\phantom{\|x_{\textrm{series}}\|_{\mathcal{H}_1}^2}
=\tr\left(\left(\sum\limits_{n=1}^N \frac{1}{a^{(n)}(K)}g_n(K)(U^*)^n\right)x^*\right)\leq
\left|\left|\sum\limits_{n=1}^N \frac{1}{a^{(n)}(K)}g_n(K)(U^*)^n\right|\right|_1||x||,
\end{gather*}
where $||x||^2_1=\tr(x^*x)$ is the trace class norm.
Next we estimate
\begin{gather*}
\left\|\sum\limits_{n=1}^N \frac{1}{a^{(n)}(K)}g_n(K)\left(U^*\right)^n\right\|_1^2 =
\tr\left(\sum\limits_{n=1}^N\frac{1}{a^{(n)}(K)}g_n(K)\left(U^*\right)^n\sum\limits_{l=1}^N
U^lg_l(K)\frac{1}{a^{(n)}(K)}\right)
\\
\qquad
=\tr\left(\sum\limits_{n=1}^N\frac{1}{a^{(n)}(K)}g_n(K)\left(U^*\right)^nU^n g_n(K)\frac{1}{a^{(n)}(K)}\right)
= \sum\limits_{k=0}^\infty \sum\limits_{n=1}^N \frac{1}{a^{(n)}(k)a^{(n)}(k)}|g_n(k)|^2
\\
\qquad
\le \|x_{\textrm{series}}\|^2_{\mathcal{H}_1}\sup_{n,k}\left(\frac{1}{a^{(n)}(k)}\right).
\end{gather*}
Using the the summability conditions on the weights we obtain
\begin{gather*}
\sup_{n,k}\left(\frac{1}{a^{(n)}(k)}\right) \le \sup_{n}\left(\sum\limits_{k=0}^\infty \frac{1}{a^{(n)}(k)}\right)=
\sup_{n} (s(n)) \le \textrm{const},
\end{gather*}
and hence we get $\|x_{\textrm{series}}\|_{\mathcal{H}_1}^2 \le
\textrm{const}\|x_{\textrm{series}}\|_{\mathcal{H}_1}\|x\|$.
This shows the continuity of the map $C^*(U)\ni x\to x_\textrm{series}\in \mathcal{H}_1$, and consequently
$x_\textrm{series}$ converges in $\mathcal{H}_1$ for every $x\in C^*(U)$.

Next we show that the map $C^*(U)\ni x\to x_\textrm{series}\in \mathcal{H}_1$ is one-to-one.
Let $x$ and $y$ belong to $C^*(U)$ and suppose that $x_\textrm{series} = y_\textrm{series}$.
This means that $x_n^+(k) = y_n^+(k)$ for all $n\ge0$ and all $k$, and $x_n^-(k) = y_n^-(k)$ for all $n\ge1$ and all
$k$.
From $x_n^+(k) = \langle e_k, (U^*)^nx e_k\rangle$ and $x_n^-(k) = \langle e_k, xU^n e_k\rangle$ it follows that all
matrix coef\/f\/icients of $x$ and $y$ are the same so we must have $x=y$.
Thus the map $C^*(U)\ni x\to x_\textrm{series}\in \mathcal{H}_1$ is one-to-one.

To prove density we def\/ine the following indicator function
\begin{gather*}
\delta_l(k) =
\begin{cases}
1, & l=k,
\\
0, & l\neq k.
\end{cases}
\end{gather*}
It is clear that $U^n\delta_l(K)$ and $\delta_l(K)(U^*)^n$ are in $C^*(U)$, and moreover they form an orthogonal basis
for $\mathcal{H}_1$.
Finally f\/inite linear combinations of $U^n\delta_l(K)$ and $\delta_l(K)(U^*)^n$ form a~dense subspace of $C^*(U)$, hence
making it a~dense subspace of $\mathcal{H}_1$.
Thus the proof is complete.
\end{proof}

We can now start the analysis of the Dirac type operators on the quantum sphere.

\begin{defin}
The mirror quantum sphere of~\cite{HMS}, denoted $C(S_{m}^2)$, is def\/ined to be the following $C^*$-algebra
\begin{gather*}
C(S_{m}^2) \cong \left\{(x,y)\in C^*(U)\times C^*(U): \sigma(x) = \overline{\sigma}(y)\right\}.
\end{gather*}
\end{defin}

The Hilbert space of the square integrable functions on the mirror quantum sphere is def\/ined to be $\mathcal{H}:=
\mathcal{H}_1\oplus \mathcal{H}_1$.

To describe the class of the Dirac type operators that we will be working with we need the following Jacobi type
operators
\begin{gather*}
\overline{A}^{(n)}f(k) = b^{(n+1)}(k)\Big(f(k) - c_+^{(n)}(k)f(k+1)\Big): \ell_{a^{(n)}}^2({\mathbb N})\to
\ell_{a^{(n+1)}}^2({\mathbb N}),
\\
A^{(n)}f(k) = b^{(n)}(k)\Big(f(k) - c_-^{(n)}(k-1)f(k-1)\Big): \ell_{a^{(n+1)}}^2({\mathbb N})\to
\ell_{a^{(n)}}^2({\mathbb N}),
\end{gather*}
where $\dom(\overline{A}^{(n)}) = \{f\in \ell_{a^{(n)}}^2({\mathbb N}):
\|\overline{A}^{(n)}f\|_{a^{(n+1)}}<\infty\}$, $\dom(A^{(n)}) = \{f\in \ell_{a^{(n+1)}}^2({\mathbb N}):
\|A^{(n)}f\|_{a^{(n)}}<\infty\}$, and where $\ell_{a^{(n)}}^2({\mathbb N})$ is the space of sequences, $\{f(k)\}$,
$k\in{\mathbb N}$, such that $\sum\limits_k \frac{1}{a^{(n)}(k)}|f(k)|^2$ is f\/inite.
(In the above formula we assume $f(-1)=0$.)

Notice that the above operators are in fact general one-step dif\/ference operators.
For our analytical purposes we will require the following conditions on the coef\/f\/icients
(the f\/irst condition is a~repeat from above):
\begin{itemize}\itemsep=0pt
\item
$\{a^{(n)}(k)\}$ is a~sequence of positive numbers such that the sum $s(n):=\sum\limits_{k=0}^\infty
\frac{1}{a^{(n)}(k)}$ exists and $s(n)$ goes to zero as $n\to\infty$.

\item
$\{b^{(n)}(k)\}$ is a~sequence of positive numbers such that the sum $t(n):= \sum\limits_{k=0}^\infty
\frac{a^{(n)}(k)}{b^{(n)}(k)^2}$ exists and $t(n)$ is bounded in $n$.

\item
$\{c_{\pm}^{(n)}(k)\}$ are sequences of real numbers such that for all $M$, $N$, and $n$ there is a~positive number
$\kappa$ independent of $M$, $N$, and $n$ such that $\kappa \le \prod\limits_{k=M}^N c_{\pm}^{(n)}(k) \le
\frac{1}{\kappa}$.
Moreover we require that the product $\prod\limits_{k=0}^\infty \frac{1}{c_{\pm}^{(n)}(k)}$ exists for each $n$.
\end{itemize}

An example of such sequences coming from the theory of quantum groups is described at the end of the paper.

\begin{defin}
Dirac type operators on the quantum unit disk are the operators $\delta$ def\/ined for formal power series $f$ by
\begin{gather}
\label{quantDTO}
\delta f = -\sum\limits_{n\ge0} U^{n+1} \overline{A}^{(n)}f_n^+(k) + \sum\limits_{n\ge1} A^{(n-1)}f_n^-(k)(U^*)^{n-1},
\end{gather}
with the Jacobi type operators $A^{(n)}$, $\overline {A}^{(n)}$ def\/ined above.
\end{defin}

The structure of $\delta$ is chosen to mimic the d-bar operator from the previous section, and it matches the structure
of the Dirac type operators of our previous papers~\cite{CKW,KM1,KM2,KM4,KM5}.

Fixing the sequences $\{a^{(n)}(k)\}$, $\{b^{(n)}(k)\}$, and $\{c_{\pm}^{(n)}(k)\}$ we can now def\/ine the Dirac type
operator $D$ on the mirror quantum sphere, the operator that we will be studying for the remainder of this paper.
We set
\begin{gather}
\label{quant_dirac_op}
D(f,g) = (\delta f, \delta g),
\end{gather}
where $\dom(D) = \{ (f,g)\in\mathcal{H}: \|\delta f\|<\infty, \|\delta g\|<\infty,\sigma(f) =
\overline{\sigma}(g)\}$.

We remark here that with some extra ef\/fort the operator above can be interpreted as acting between dif\/ferent quantum
line bundles, by checking the coef\/f\/icients of the Fourier decomposition of $f$ and $\delta f$.
The description of the quantum line bundles on the quantum 2-sphere that uses gluing from two disks can be found
in~\cite{Wagner}, which can presumably be adapted to the mirror quantum sphere studied here.

The following proposition is a~direct consequence of formulas~\eqref{lim_at_infty1} and~\eqref{lim_at_infty2}.

\begin{prop}
\label{bndy_cond_quant}
The boundary condition $ \sigma(f) = \overline{\sigma}(g)$ is equivalent to the following: if $n\ge 1$ then
$f_n^+(\infty) = g_n^-(\infty)$,
$f_n^-(\infty) = g_n^+(\infty)$ and $f_0^+(\infty) = g_0^+(\infty)$.
\end{prop}

From this proposition we can rewrite the domain for $D$ in the following way
\begin{gather*}
\dom(D) = \{(f,g)\in\mathcal{H}: \|\delta f\|<\infty,~\|\delta g\|<\infty,~\text{such that}~(1)~\text{and}~(2)~\text{hold}\},
\end{gather*}
where
\begin{gather*}
(1)\quad\text{If}\quad n\ge 1 \quad\text{then}\quad  f_n^+(\infty) = g_n^-(\infty)\quad\text{and}\quad  f_n^-(\infty) = g_n^+(\infty),
\\
(2)\quad f_0^+(\infty) = g_0^+(\infty).
\end{gather*}
We can now state the main result of this paper.

\begin{theo}
\label{compact_para_quant}
The Dirac type operator $D$, defined in equation~\eqref{quant_dirac_op}, subject to the boundary conditions given in
Proposition~{\rm \ref{bndy_cond_quant}}, is a~right invertible operator and is left invertible modulo compact operators.
Moreover $D$ has a~compact parametrix.
\end{theo}

The theorem is proved in a~sequence of steps.
First we study the kernel of the Jacobi opera\-tors~${A}^{(n)}$ and~$\overline{A}^{(n)}$ that appear in the def\/inition of~$D$.
\begin{lem}
\label{A_and_A_bar_ker}
The kernel of $A^{(n)}$ is trivial and for some constant $\alpha$ we have
\begin{gather*}
\Ker \overline{A}^{(n)} = \left\{\left(\prod\limits_{i=0}^{k-1}\frac{1}{c_+^{(n)}(i)}\right)\alpha\right\}.
\end{gather*}
\end{lem}

\begin{proof}
The goal is solve the following two equations $A^{(n)}f(k) = 0$ and $\overline{A}^{(n)}f(k)=0$.
For the f\/irst equation, since we take $f(-1)=0$, it is easy to see inductively that $f(k)=0$ for every $k$, so that
$A^{(n)}$ has no kernel.

Using the formula for $\overline{A}^{(n)}$ and solving recursively it is easy to see that $\overline{A}^{(n)}f(k)=0$
produces
\begin{gather*}
f(k) = \left(\prod\limits_{i=0}^{k-1}\frac{1}{c_+^{(n)}(i)}\right)f(0).
\end{gather*}
We must show ${f(k)}\in \ell_{a^{(n)}}^2({\mathbb N})$.
A simple computation yields
\begin{gather*}
\|f\|_{a^{(n)}}^2 = \sum\limits_{k=0}^{\infty}
\frac{1}{a^{(n)}(k)}\left(\prod\limits_{i=0}^{k-1}\frac{1}{c_+^{(n)}(i)}\right)^2|f(0)|^2
\leq\frac{s(n)|f(0)|^2}{\kappa^2} < \infty
\end{gather*}
because of the conditions on $c_+^{(n)}(k)$ and $a^{(n)}(k)$.
Thus the proof is complete.
\end{proof}
Next we look at the nonhomogeneous equations for $A^{(n)}$ and $\overline{A}^{(n)}$.

\begin{lem}
\label{A_and_A_bar_gen_soln}
The solution to the equation $A^{(n)}f(k) = g(k)$ is
\begin{gather*}
f(k) = \sum\limits_{i=0}^k\frac{1}{b^{(n)}(i)}\left(\prod\limits_{j=i}^{k-1}c_-^{(n)}(j)\right)g(i).
\end{gather*}
Moreover the solutions to the equation $\overline{A}^{(n)}f(k) = -g(k)$ are
\begin{gather*}
f(k) = \left(\prod\limits_{i=0}^{k-1}\frac{1}{c_+^{(n)}(i)}\right)\alpha - \sum\limits_{i=k}^\infty
\frac{1}{b^{(n+1)}(i)}\left(\prod\limits_{j=k}^{i-1}c_+^{(n)}(j)\right)g(i),
\end{gather*}
where $\alpha$ is an arbitrary constant.
\end{lem}

\begin{proof}
The formulas follow from simple calculations, see the proofs in Proposition~4.2 and Proposition~4.4 in~\cite{KM1}.
We need to show that $f(k)$ in each case belongs to its respective $\ell^2$ space.
For the solution $f(k)$ of $A^{(n)}f(k) = g(k)$ we estimate pointwise
\begin{gather*}
|f(k)| \le \frac{1}{\kappa}\sum\limits_{i=0}^k\frac{\sqrt{a^{(n)}(i)}}{b^{(n)}(i)}\frac{g(i)}{\sqrt{a^{(n)}(i)}}
\end{gather*}
by using the condition on $c_-^{(n)}(k)$.
The Cauchy--Schwarz inequality implies
\begin{gather*}
|f(k)|^2 \le \frac{1}{\kappa^2}\sum\limits_{i=0}^k \frac{a^{(n)}(i)}{b^{(n)}(i)^2}\sum\limits_{i=0}^k
\frac{|g(i)|^2}{a^{(n)}(i)}\le \frac{t(n)}{\kappa^2}\|g\|_{a^{(n+1)}}^2,
\end{gather*}
and thus we have
\begin{gather*}
\|f\|_{a^{(n)}}^2 = \sum\limits_{k=0}^\infty \frac{1}{a^{(n)}(k)}|f(k)|^2 \le
\frac{t(n)}{\kappa^2}\|g\|_{a^{(n+1)}}^2\sum\limits_{k=0}^\infty \frac{1}{a^{(n)}(k)}=
\frac{t(n)s(n)}{\kappa^2}\|g\|_{a^{(n+1)}}^2 < \infty
\end{gather*}
by the assumptions on $s(n)$ and $t(n)$.
A similar argument shows that the $f(k)$ solving $\overline{A}^{(n)}f(k) = -g(k)$ is in $\ell_{a^{(n+1)}}^2({\mathbb
N})$.
Thus this completes the proof.
\end{proof}

The above lemma allows us to study the limits at inf\/inity.

\begin{lem}
Let $f\in \dom(A^{(n)})$, then $f(\infty):= \lim\limits_{k\to\infty} f(k)$ exists.
Moreover if $f\in \dom(\overline{A}^{(n)})$, then $f(\infty):= \lim\limits_{k\to\infty} f(k)$ exists.
\end{lem}

\begin{proof}
If $f\in \dom(A^{(n)})$, then by using Lemma~\ref{A_and_A_bar_gen_soln} we have
\begin{gather*}
f(\infty) = \sum\limits_{i=0}^\infty
\frac{1}{b^{(n)}(i)}\left(\prod\limits_{j=i}^{\infty}c_-^{(n)}(j)\right)A^{(n)}f(i).
\end{gather*}
Similarly, if $f\in\dom\big(\overline{A}^{(n)}\big)$ then $f(\infty) = \left(\prod\limits_{i=0}^\infty 1/c_+^{(n)}(i)
\right)\alpha$.
Because of the conditions on $b^{(n)}(k)$ and $c_-^{(n)}(k)$ it is clear that the inf\/inite products and sums in the
above formulas are convergent.
Thus the proof is complete.
\end{proof}

It will be advantageous to write the solution to $\overline{A}^{(n)}f(k) = -g(k)$ in a~slightly dif\/ferent way.
This will not only make the estimates simpler, it will also mimic the continuous case theory very closely.
Using the above we get
\begin{gather*}
f(k) = \left(\prod\limits_{i=k}^\infty c_+^{(n)}(i)\right)f(\infty) - \sum\limits_{i=k}^\infty
\frac{1}{b^{(n+1)}(i)}\left(\prod\limits_{j=k}^{i-1}c_+^{(n)}(j)\right)g(i).
\end{gather*}

Now we look at the kernel of $D$.

\begin{prop}
The kernel of $D$ subject to the boundary condition given in Proposition~{\rm \ref{bndy_cond_quant}}, is one-dimensional.
\end{prop}

\begin{proof}
We wish to study the equation $D(x,y) = 0$, with the Fourier decomposition for $x= \sum\limits_{n\ge0} U^n x_n^+(K) +
\sum\limits_{n\ge1} x_n^-(K)(U^*)^n$, and similarly for y.
This translates to the following system of equations
\begin{gather*}
\delta x = 0,
\qquad
\delta y = 0
\end{gather*}
subject to the boundary condition.
Using the Fourier decomposition and Lemma~\ref{A_and_A_bar_ker}, it follows that $x_n^-(k) = 0$ for all $k$ and
\begin{gather*}
x_n^+(k) = \left(\prod\limits_{i=k}^\infty c_+^{(n)}(i)\right)x_n^+(\infty).
\end{gather*}
We also have identical formulas for $y$.
Since $x_n^-(k) = 0$ and $y_n^-(k) = 0$ for all $k$ it follows that $x_n^-(\infty)=0$ and $y_n^-(\infty)= 0$.
If $n\ge1$ then from the boundary condition we get $0 = y_n^-(\infty) = x_n^+(\infty)$ which implies that $x_n^+(k)=0$
for all $k$.
The same type of argument shows that $y_n^+(k)=0$ for all $k$.
If $n=0$ we get
\begin{gather*}
x_0^+(k) = \left(\prod\limits_{i=k}^\infty c_+^{(n)}(i)\right)x_0^+(\infty)
\qquad
\textrm{and}
\qquad
y_0^+(k) = \left(\prod\limits_{i=k}^\infty c_+^{(n)}(i)\right)y_0^+(\infty).
\end{gather*}
However we have $x_0^+(\infty) = y_0^+(\infty)$ which implies in particular, that $x_0^+(0) = y_0^+(0)$ is an arbitrary
constant.
Therefore it follows that $\Ker D$ is one-dimensional.
\end{proof}

The next goal is to construct the parametrix of the Dirac type operator subject to the boundary condition.
First we def\/ine three integral operators that the parametrix will decompose into in a~similar way to the continuous
case.
The operators are
\begin{gather}
T_1^{(n)}f(k) = \left(\prod\limits_{i=k}^\infty c_+^{(n)}(i)\right)\sum\limits_{i=0}^\infty
\frac{a^{(n-1)}(i)}{b^{(n-1)}(i)}\left(\prod\limits_{j=i}^\infty c_-^{(n-1)}(j)\right)\frac{f(i)}{a^{(n-1)}(i)},
\nonumber
\\
T_2^{(n)}f(k) = -\sum\limits_{i=0}^\infty\frac{a^{(n+1)}(i)}{b^{(n+1)}(i)}\chi\left(\frac{k}{i}\right)
\left(\prod\limits_{j=k}^{i-1}c_+^{(n)}(j)\right)\frac{f(i)}{a^{(n+1)}(i)},
\nonumber
\\
T_3^{(n)}f(k) = \sum\limits_{i=0}^\infty
\frac{a^{(n-1)}(i)}{b^{(n-1)}(i)}\chi\left(\frac{i}{k}\right)\left(\prod\limits_{j=i}^{k-1}c_-^{(n-1)}(j)\right)\frac{f(i)}{a^{(n-1)}(i)},
\label{int_ops_quant}
\end{gather}
where $T_1^{(n)}, T_3^{(n)}: \ell_{a^{(n-1)}}^2({\mathbb N})\to\ell_{a^{(n)}}^2({\mathbb N})$ and $T_2^{(n)}:
\ell_{a^{(n+1)}}^2({\mathbb N}) \to \ell_{a^{(n)}}^2({\mathbb N})$ and again $\chi(t)=1$ for $t\le1$ and zero otherwise.

We will also need the following rank one operator in $\mathcal{H}$ def\/ined by $C(x,y) = (\tilde x(K),\tilde x(K))$,
where
\begin{gather*}
\tilde x(k) = \left(\prod\limits_{i=0}^{k-1}\frac{1}{c_+^{(n)}(i)}\right)x_0^+(0)
\end{gather*}
with $x= \sum\limits_{n\ge0} U^n x_n^+(K) + \sum\limits_{n\ge1} x_n^-(K)(U^*)^n$.

\begin{prop}
\label{quant_para}
The Dirac type operator $D$, defined in equation~\eqref{quant_dirac_op}, subject to the boundary condition given in
Proposition~{\rm \ref{bndy_cond_quant}}, is right invertible, with the right inverse $Q$ given by $Q(p,q)=({x},{y})$ where
\begin{gather*}
{x} = T_2^{(0)}p_1^+(k) + \sum\limits_{n\ge 1}U^n\left(T_2^{(n)}p_{n+1}^+(k) + T_1^{(n)}q_{n-1}^-(k)\right) +
\sum\limits_{n\ge1} T_3^{(n)}p_{n-1}^-(k)(U^*)^n,
\\
{y} = T_2^{(0)}q_1^+(k) + \sum\limits_{n\ge 1}U^n\left(T_2^{(n)}q_{n+1}^+(k) + T_1^{(n)}p_{n-1}^-(k)\right) +
\sum\limits_{n\ge1} T_3^{(n)}q_{n-1}^-(k)(U^*)^n.
\end{gather*}
Moreover we have $QD = I - C$.
\end{prop}

\begin{proof}
The goal is to solve the equation $(\delta x, \delta y) = (p,q)$.
Using the Fourier decompositions, the equation $\delta x = p$ will produce two type of equations:
$\overline{A}^{(n)}x_n^+(k) = -p_{n+1}^+(k)$ and $A^{(n-1)}x_n^-(k) = q_{n-1}^-(k)$.
Using Lemma~\ref{A_and_A_bar_gen_soln}, we immediately get
\begin{gather*}
x_n^+(k) = \left(\prod\limits_{i=k}^\infty c_+^{(n)}(i)\right)x_n^+(\infty) -
\sum\limits_{i=k}^\infty\frac{1}{b^{(n+1)}(i)}\left(\prod\limits_{j=k}^{i-1}c_+^{(n)}(j)\right)p_{n+1}^+(i),
\\
x_n^-(k) = \sum\limits_{i=0}^k \frac{1}{b^{(n-1)}(i)}\left(\prod\limits_{j=i}^{k-1} c_-^{(n-1)}(j)\right)p_{n-1}^-(i),
\end{gather*}
and the formulas for $y_n^+(k)$ and $y_n^-(k)$ are the same except we need to replace the $p$'s with $q$'s.
Next we apply the boundary condition for $n\ge 1$.
Using $x_n^+(\infty) = y_n^-(\infty)$ we get
\begin{gather*}
x_n^+(\infty) = \sum\limits_{i=0}^\infty \frac{1}{b^{(n-1)}(i)}\left(\prod\limits_{j=i}^\infty
c_-^{(n-1)}(j)\right)q_{n-1}^-(i).
\end{gather*}
Similarly applying $y_n^+(\infty) = x_n^-(\infty)$ we obtain
\begin{gather*}
y_n^+(\infty) = \sum\limits_{i=0}^\infty \frac{1}{b^{(n-1)}(i)}\left(\prod\limits_{j=i}^\infty
c_-^{(n-1)}(j)\right)p_{n-1}^-(i).
\end{gather*}
Using the integral formulas~\eqref{int_ops_quant}, it's not too hard to see that we get
\begin{gather*}
x_n^+(k) = T_1^{(n)}q_{n-1}^-(k) + T_2^{(n)}p_{n+1}^+(k),
\qquad
x_n^-(k) = T_3^{(n)}p_{n-1}^-(k),
\\
y_n^+(k) = T_1^{(n)}p_{n-1}^-(k) + T_2^{(n)}q_{n+1}^+(k),
\qquad
y_n^-(k) = T_3^{(n)}q_{n-1}^-(k).
\end{gather*}
The above formulas for the coef\/f\/icients imply the formulas for ${x}$ and ${y}$.
From the construction we have $DQ(p,q)=(p,q)$.
As with the continuous case, a~direct substitution using the boundary conditions shows that $QD(x,y) = (I-C)(x,y)$.
Thus the proof is complete.
\end{proof}

It follows immediately from the above proof that the $C$ def\/ined in the above proposition is a~f\/inite rank operator
making $D$ left invertible modulo compacts.
The next proposition is the key step to show that the $Q$ we just constructed in the quantum case is a~compact operator.

\begin{prop}
\label{hs_norm_quant_T}
The three integral operators given in equation~\eqref{int_ops_quant} are Hilbert--Schmidt operators.
Moreover the Hilbert--Schmidt norms of these operators go to zero as $n\to\infty$.
\end{prop}

\begin{proof}
We f\/irst look at the Hilbert--Schmidt norm of $T_3^{(n)}$.
We have
\begin{gather*}
\|T_3^{(n)}\|_{\textrm{HS}}^2 = \sum\limits_{k=0}^\infty \sum\limits_{i=0}^\infty
\frac{1}{a^{(n)}(k)}\frac{1}{a^{(n-1)}(i)}\left(\frac{a^{(n-1)}(i)}{b^{(n-1)}(i)}\prod\limits_{j=i}^{k-1}c_-^{(n-1)}(j)\right)^2
\chi\left(\frac{i}{k}\right)
\\
\phantom{\|T_3^{(n)}\|_{\textrm{HS}}^2}
\le \frac{1}{\kappa^2}\sum\limits_{k=0}^\infty
\frac{1}{a^{(n)}(k)}\sum\limits_{i=0}^\infty\frac{a^{(n-1)}(i)}{b^{(n-1)}(i)^2} = \frac{s(n)t(n-1)}{\kappa^2}.
\end{gather*}

The same arguments work for $T_1^{(n)}$ and $T_2^{(n)}$ leading to
\begin{gather*}
\|T_1^{(n)}\|_{\textrm{HS}}\le \frac{\sqrt{s(n)t(n-1)}}{\kappa^2}
\qquad
\textrm{and}
\qquad
\|T_2^{(n)}\|_{\textrm{HS}}\le \frac{\sqrt{s(n)t(n+1)}}{\kappa}.
\end{gather*}
Since $s(n)$ goes to zero as $n\to\infty$ and $t(n)$ is bounded in $n$, it follows that all the Hilbert--Schmidt norms
go to zero as $n\to\infty$.
Thus this completes the proof.
\end{proof}

\begin{proof}[Proof of Theorem~\ref{compact_para_quant}] We repeat the steps of the proof of Theorem~\ref{compact_para_cont}.
Namely Proposition~\ref{quant_para} shows that~$Q$ is the right inverse to~$D$.
Moreover the comment directly after the proof of Proposition~\ref{quant_para} shows that~$D$ is left invertible modulo
compacts.
From the decompostion of~$Q$ in Proposition~\ref{quant_para}, and from the fact that the Hilbert--Schmidt norms of~$T_1^{(n)}$, $T_2^{(n)}$, and~$T_3^{(n)}$ go zero as $n\to\infty$ by Proposition~\ref{hs_norm_quant_T}, it follows that~$Q$ must be a~compact operator because it is an inf\/inite direct sum of compact operators with norms going to zero.
This ends the proof.
\end{proof}

We end this paper with an example of a~Dirac type operator on the quantum sphere obtained by gluing the quantum d-bar
operator of~\cite{KM1} with weights of~\cite{KL1}.
Let $0\le q <1$.
Def\/ine the $q$-weight $w(k)$ by $w^2(k) = 1 - q^{k+1}$, see~\cite{KL1}.
It is clear that these weights monotonically increase to $1$.
Next let $U_W$ be the weighted shift with weight $w(k)$, in other words, $U_We_k = w(k)e_{k+1}$, and def\/ine $S(k) =
w^2(k) - w^2(k-1)$, so that $S(K)=[U_W^*,U_W]$.
It is clear that we have $\sum\limits_{k=0}^\infty S(k) = 1$.
In our example the operator $\delta$ on the quantum disk is given by a~commutation with the above weighted unilateral
shift: $\delta x = S(K)^{-1/2}[x,U_W]S(K)^{-1/2}$, see~\cite{KM1}.
Then we have the following decomposition of $\delta$
\begin{gather*}
\delta x = -\sum\limits_{n=0}^\infty U^{n+1}\overline{A}^{(n)}x^+_n(k) + \sum\limits_{n=1}^\infty
A^{(n-1)}x^-_n(k)(U^*)^{n-1},
\end{gather*}
where again $x = \sum\limits_{n=0}^\infty U^n x^+_n(k) + \sum\limits_{n=1}^\infty x^-_n(k)(U^*)^n$.
The Hilbert space norm is given by the expression
$\|x\|_{\mathcal{H}_1}^2={{\operatorname{Tr}}}\left(S(K)^{1/2}xS(K)^{1/2}x^*\right)$.

A direct calculation shows that $c_+^{(n)}(k) = w(k)/w(k+n)$, $c_-^{(n)}(k) = (w(k)w(k+n))/w^2(k+n+1)$, $a^{(n)}(k) =
S^{-1/2}(k)S^{-1/2}(k+n)$, and $b^{(n)}(k) = a^{(n-1)}(k)w(k+n-1)$.
It is clear that $s(n)$ is f\/inite and $s(n) = \sum\limits_{k=0}^\infty \frac{1}{a^{(n)}(k)} =
q^{n/2}(1-q)\sum\limits_{k=0}^\infty q^k = q^{n/2}$.
This calculation also shows that $s(n)\to 0$ as $n\to\infty$.

Next we examine $t(n)$.
We have
\begin{gather*}
t(n) = \sum\limits_{k=0}^\infty \frac{a^{(n)}(k)}{b^{(n)}(k)^2} = \frac{(1-q)q^{n/2}}{q}\sum\limits_{k=0}^\infty
\frac{q^k}{1-q^{k+n}} \le q^{(n-2)/2}\sum\limits_{k=0}^\infty q^k = \frac{q^{(n-2)/2}}{1-q}
\end{gather*}
which is bounded in $n$.

Since $w(k)$ is monotonically increasing, it follows that $c_{\pm}^{(n)}(k)\le 1$ for all $k$ and $n$.
We compute
\begin{gather*}
\prod\limits_{k=0}^\infty \frac{1}{c_+^{(n)}(k)} = \prod\limits_{k=0}^\infty \sqrt{\frac{1-q^{k+n+1}}{1-q^{k+1}}} =
\frac{1}{\sqrt{(1-q^2)\cdots(1-q^{n+1})}} < \infty
\end{gather*}
for all $n$.
Similarly we get
\begin{gather*}
\prod\limits_{k=0}^\infty \frac{1}{c_-^{(n)}(k)} = \prod\limits_{k=0}^\infty \sqrt{\frac{1-q^{k+n+2}}{1-q^{k+1}}}
\sqrt{\frac{1-q^{k+n+2}}{1-q^{k+n+1}}} = \frac{1}{\sqrt{(1-q)(1-q^2)\cdots(1-q^{n+1})}} < \infty
\end{gather*}
for all $n$.
Therefore for $\kappa=1-q$ we have
\begin{gather*}
\kappa \le \prod\limits_{k=M}^N c_{\pm}^{(n)}(k) \le \frac{1}{\kappa}
\end{gather*}
for all $M$, $N$, and $n$.
As a~consequence all the coef\/f\/icient conditions are satisf\/ied and the operator $\delta$ yields a~well behaved Dirac type
operator on the mirror quantum sphere.

\subsection*{Acknowledgements}

The f\/irst author would like to thank W.~Szymanski for interesting discussions that inspired this paper.
We would also like to acknowledge the referees for insightful comments and remarks.

\pdfbookmark[1]{References}{ref}
\LastPageEnding

\end{document}